\documentclass[11pt]{article} 

\usepackage{makeidx}         
\usepackage{graphicx}        
\usepackage{multicol}        
\usepackage[bottom]{footmisc}

\usepackage{amssymb, amsmath}
\usepackage{amsthm}
\usepackage{mathtools} 

\usepackage{authblk}
\usepackage[margin=0.7in]{geometry}

\newcommand{\N}{\mathbb{N}}
\newcommand{\Z}{\mathbb{Z}}
\newcommand{\R}{\mathbb{R}}
\newcommand{\Sum}{\displaystyle\sum}

\newcommand{\KB}{\mathcal{S}_{2,\infty}}
\newcommand{\RR}{\R P^3\sharp\R P^3}
\newcommand{\RP}{\R P}
\newcommand{\Om}{\Omega}

\newcommand{\calS}{\mathcal{S}}
\newcommand{\calL}{\mathcal{L}}

\newcommand{\calB}{\mathcal{B}}

\newcommand{\kbsm}{\mathcal{S}_{2,\infty}}
\newcommand{\hsm}{\mathcal{S}_3}

\usepackage{extarrows}

\usepackage[percent]{overpic} 
\usepackage[caption=false]{subfig} 
\usepackage{enumerate} 

\newtheorem{theorem}{Theorem}
\newtheorem{proposition}{Proposition}
\newtheorem{lemma}{Lemma}
\newtheorem{example}{Example}
\newtheorem{definition}{Definition}

\makeindex

\begin{document}

\date{}

\title{Link diagrams in Seifert manifolds and applications to skein modules}

\author[1]{Bo\v stjan Gabrov\v sek}
\author[2]{Maciej Mroczkowski}

\affil[1]{\normalsize Faculty of Mathematics and Physics, University of Ljubljana, Jadranska 21, 1000 Ljubljana, Slovenia.}
\affil[2]{\normalsize Institute of Mathematics, Faculty of Mathematics, Physics and Informatics, University of Gdansk, 80-308 Gdansk, Poland.}




\maketitle


\abstract{In this survey paper we present results about link diagrams in Seifert
manifolds using arrow diagrams, starting with link diagrams in $F\times S^1$ 
and $N\hat{\times}S^1$, where $F$ is an orientable and $N$ an unorientable
surface. 
Reidemeister moves for such arrow diagrams make the study of link
invariants possible. Transitions between arrow diagrams and alternative diagrams are presented. We recall results about 
the Kauffman bracket and HOMFLYPT skein modules of some Seifert manifolds using arrow diagrams, namely lens spaces, 
a product of a disk with two holes times $S^1$, $\R P^3 \# \R P^3$, and prism manifolds. 
We also present new bases of the Kauffman bracket and HOMFLYPT skein modules of the solid torus and lens spaces.}

\section{Arrow diagrams of links in products and twisted products of $S^1$ 
and a surface}
\label{sec:arrow_diagrams}

Let $F$ be an orientable surface and $N$ an unorientable surface.
In this section we recall the construction of arrow diagrams for links
in $F\times S^1$, introduced in \cite{MD},
 and $N\hat{\times} S^1$, introducted in \cite{M1}. These diagrams
are very similar to gleams introduced in \cite{T2}.

\subsection{Arrow diagrams of links in $F\times S^1$}
Let $L$ be a link in $M=F\times S^1$. We cut $M$ along $F_0=F\times\{1\}$,
$1\in S^1$, to get $M'=F\times [0,1]$. By a general position argument we may
assume that $L$ intersects $F_0$ transversally in a finite number of points.
In $M'$ the link $L$ becomes $L'$ - a collection of circles and arcs with 
endpoints coming in pairs $(x,0)$ and $(x,1)$, $x\in F$.
Let $\pi$ be the vertical projection from $M'$ onto $F$.
Then  $\pi(L')$ is a collection of closed curves. Again, by a general position 
argument, we may assume that there are only transversal double points
in $\pi(L')$ and the endpoints of arcs are projected onto points distinct
from these double points. An arrow diagram $D$ of the link $L$ is $\pi(L')$
with some extra information: for double points $P$, the usual information of
over- and undercrossing is encoded depending on the relative height of the
two points $\pi^{-1}(P)$ in $F\times[0,1]$; for points $Q$ that are projections
of endpoints of arcs $(x,0)$ and $(x,1)$ in $L'$, orient $L$ in such a way
that the height drops by $1$ in $L'$ when the first coordinate crosses $x$, 
and put on $Q$ an arrow indicating this orientation. 

\begin{figure}[h]
	\centering\hspace{0.8cm}
	\begin{overpic}[page=1]{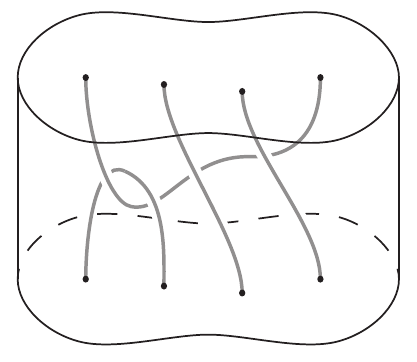}
		\put(-23,15){$F\times\{0\}$}
		\put(-23,65){$F\times\{1\}$}
	\end{overpic}\raisebox{5.5em}{$\xrightarrow{\;\;\pi\;\;}$}
	\raisebox{1.5em}{\begin{overpic}[page=2]{images}
	\end{overpic}}
	\caption{A link in $F\times S^1$ and its diagram.}\label{fig:diagram}
\end{figure}

Thus, an arrow diagram $D$ is a collection of immersed curves in $F$, with
under- and overcrossing information for double points and some arrows on these
curves. 
For an example see the diagram on Figure~\ref{fig:diagram}.

\begin{figure}[h]
	\centering
	\subfloat[]{\centering
	\includegraphics[page=46]{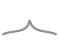}\label{fig:forbid1}
	}\hspace{0.8cm}
	\subfloat[]{\centering
	\includegraphics[page=47]{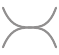}\label{fig:forbid2}
	}\hspace{0.8cm}
	\subfloat[]{\centering
	\includegraphics[page=48]{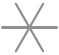}\label{fig:forbid3}
	}\hspace{0.8cm}
	\subfloat[]{\centering
	\includegraphics[page=54]{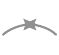}\label{fig:forbid4}
	}\hspace{0.8cm}
	\subfloat[]{\centering
	\includegraphics[page=49]{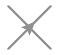}\label{fig:forbid5}
	}
	\caption{Forbidden positions of regular diagrams in $F\times S^1$.}
	\label{fig:forbidden1}
\end{figure}

We call an arrow diagram regular, if none of the following forbidden positions appear on the diagram:

\begin{enumerate}[i)]
	\item cusps~(Fig.~\ref{fig:forbid1}),\hfill {(1)}
	\item self-tangency points~(Fig.~\ref{fig:forbid2}),
	\item triple points~(Fig.~\ref{fig:forbid3}),
	\item two arrows coincide~(Fig.~\ref{fig:forbid4}),
	\item arrows and crossings coincide~(Fig.~\ref{fig:forbid5}).
\end{enumerate}

With standard arguments of general position we may assume that every link admits a regular diagram.

We complete this section by providing a list of generalized Reidemeister moves associated with the arrow diagrams. 
As usual, Reidemeister moves coincide with the change of diagram that occurs when an isotopy of the link is performed in such a way that the configuration of arcs (and arrows) passes through a forbidden position in the projection. 
We therefore assign each forbidden position an associated Reidemeister move. By general position theory, 
the ambient isotopy bringing one link diagram to another one passes through only a finite number of such forbidden positions, hence a finite number of Reidemeister-type moves. 

\begin{figure}[h]
	\centering
	\subfloat[$\Om_1$]{\centering
	\includegraphics[page=9]{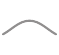}\raisebox{10pt}{$\;\;\longleftrightarrow\;\;$}\includegraphics[page=10]{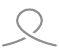}
	\label{fig:reid1}
	}\hspace{0.8cm}
	\subfloat[$\Om_2$]{\centering
	\includegraphics[page=11]{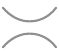}\raisebox{10pt}{$\;\;\longleftrightarrow\;\;$}\includegraphics[page=12]{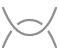}
	\label{fig:reid2}
	}\hspace{0.8cm}
	\subfloat[$\Om_3$]{\centering
	\includegraphics[page=13]{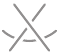}\raisebox{10pt}{$\;\;\longleftrightarrow\;\;$}\includegraphics[page=14]{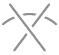}
	\label{fig:reid3}
	}\hspace{0.8cm}
	\subfloat[$\Om_4$]{\centering
	\includegraphics[page=15]{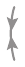}\raisebox{12pt}{$\;\;\longleftrightarrow\;\;$}\includegraphics[page=16]{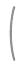}\raisebox{12pt}{$\;\;\longleftrightarrow\;\;$}\includegraphics[page=17]{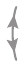}
	\label{fig:reid4}
	}\hspace{0.8cm}
	\subfloat[$\Om_5$]{\centering
	\includegraphics[page=18]{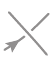}\raisebox{12pt}{$\;\;\longleftrightarrow\;\;$}\includegraphics[page=19]{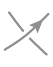}
	\label{fig:reid5}
	}
	\caption{Classical Reidemesiter moves $\Om_1$ -- $\Om_3$ and two ``arrow'' moves $\Om_4$ and $\Om_5$.}
	\label{fig:reidemeister1}
\end{figure}

\begin{figure}[t]
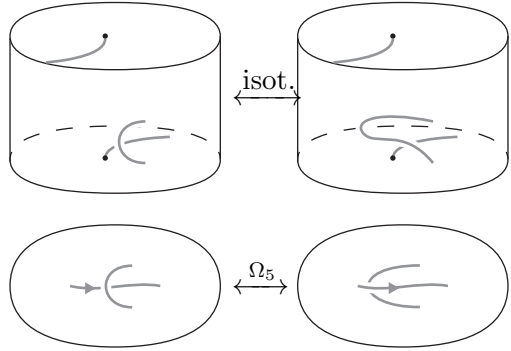

	\centering
    \begin{overpic}[page=25]{images}
   \put(45,48){$\xleftrightarrow{\mbox{isot.}}$}
   \put(45,13){$\xleftrightarrow{\;\Om_5\;}$}  
    \end{overpic}
\caption{Interpretation of the move $\Om_5$.}\label{fig:omega5}
\end{figure}

As in the classical case, Reidemeister moves $\Om_1$, $\Om_2$, and $\Om_3$ (Figure~\ref{fig:reidemeister1}) arise from the forbidden positions i), ii), and iii), respectively.
Positions iv) and v) generate Reidemeister moves $\Om_4$ (``arrow annihilation'') and $\Om_5$ (``arrow push''), respectively.
Graphical interpretation of moves $\Om_4$ and $\Om_5$ are presented in Figure~\ref{fig:omega4} and Figure~\ref{fig:omega5}, respectively.

\begin{figure}[h]
	\centering
    \begin{overpic}[page=24]{images}
   \put(29,32){$\xleftrightarrow{\mbox{isot.}}$}
   \put(64,32){$\xleftrightarrow{\mbox{isot.}}$}
    \put(29.5,9){$\xleftrightarrow{\;\Om_4\;}$}
   \put(64,9){$\xleftrightarrow{\;\Om_4\;}$}  
    \end{overpic}
\caption{Interpretation of the move $\Om_4$.}\label{fig:omega4}
\end{figure}

We conclude this section with the following Reidemeister-type theorem.

\begin{theorem}\label{thm:reid}
A link $L_1$ is ambient isotopic to a link $L_2$ if and only if an arrow diagram
$D_1$ of $L_1$ can be obtained from an arrow diagram $D_2$ of $L_2$ by a finite series of Reidemeister moves $\Om_1$ to $\Om_5$.  
\end{theorem}

\subsection{Arrow diagrams of links in $N\hat{\times} S^1$}

An unorientable surface $N$ is obtained from a sphere with $n$ holes $S_n$ by
glueing $k$ of the boundary $S^1$'s with antipodal maps (which is equivalent to
glueing Mobius bands to these holes). Denote the $k$ boundary $S^1$'s by $C$.
Let $M=N\hat{\times}S^1$ be obtained from $M'=S_n\times S^1$ by glueing 
$(x,y)\in C\times S^1$ to $(-x,r(y))$, where $r$ is a reflection of $S^1$ 
(one may take the complex conjugation). Let $L$ be a link in $M$. 
In $M'$, $L$ becomes $L'$, a collection of circles and arcs 
with endpoints coming in antipodal pairs $(x,y)$ and $(-x,r(y))$ in $C\times S^1$.

\begin{figure}[t]
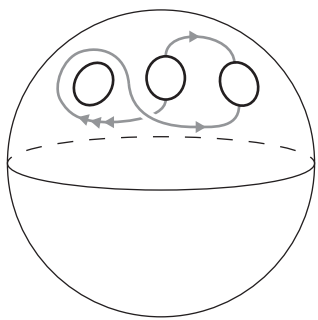

	\centering
    \begin{overpic}[page=26]{images}
    \end{overpic}
\caption{A diagram of a link in $N\hat{\times} S^1$.}\label{fig:diagram-twist}
\end{figure}

For $L'$ in $S_n\times S^1$ one constructs an arrow diagram as in the
previous subsection, the only difference being that there are now some arcs
with endpoints coming in antipodal pairs.
For an example of a diagram see Figure~\ref{fig:diagram-twist}.

\begin{figure}[h]
\centering
\subfloat[]{\centering
\includegraphics[page=47]{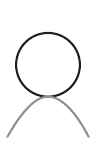}\label{fig:forbid6}
}\hspace{0.8cm}
\subfloat[]{\centering
\includegraphics[page=48]{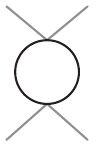}\label{fig:forbid7}
}\hspace{0.8cm}
\subfloat[]{\centering
\includegraphics[page=49]{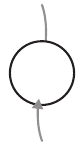}\label{fig:forbid8}
}
\caption{Additional forbidden positions of regular diagrams in $N
\hat\times S^1$.}
\label{fig:forbidden2}
\end{figure}

We call a diagram regular, if, in addition to the list (1), none of the following forbidden positions appear on the diagram:

\begin{enumerate}[i)]\setcounter{enumi}{5}
	\item tangency points with the boundary~(Fig.~\ref{fig:forbid6}),
	\item crossing coincides with the boundary~(Fig.~\ref{fig:forbid7}),
	\item arrow coincides with the boundary~(Fig.~\ref{fig:forbid8}).
\end{enumerate}

\begin{figure}[h]
	\centering
	\subfloat[$\Om_6$]{\centering
	\includegraphics[page=27]{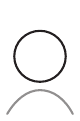}\raisebox{29pt}{$\!\longleftrightarrow\!$}\includegraphics[page=28]{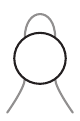}
	\label{fig:reid6}
	}\hspace{0.8cm}
	\subfloat[$\Om_7$]{\centering
	\includegraphics[page=29]{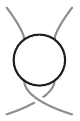}\raisebox{29pt}{$\!\longleftrightarrow\!$}\includegraphics[page=30]{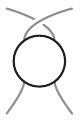}
	\label{fig:reid7}
	}\hspace{0.8cm}
	\subfloat[$\Om_8$]{\centering
	\includegraphics[page=31]{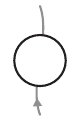}\raisebox{29pt}{$\!\longleftrightarrow\!$}\includegraphics[page=32]{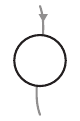}
	\label{fig:reid8}
	}
    \caption{Additional Reidemeister moves.}\label{fig:reid-twist}
\end{figure}

To these additional forbidden positions vi), vii), and viii), we associate Reidemeister moves $\Om_6$, $\Om_7$, and $\Om_8$ (Figure~\ref{fig:reid-twist}), respectively.

\begin{theorem}\label{thm:reid2}
A link $L_1$ is ambient isotopic to a link $L_2$ in $N \hat\times S^1$ if and only if an arrow diagram
$D_1$ of $L_1$ can be obtained from an arrow diagram $D_2$ of $L_2$ by a series of Reidemeister moves $\Om_1$ to $\Om_8$.  
\end{theorem}

\begin{example}\label{ex:rp3rp3}

The connected sum of two projective spaces $\RR$ is also a twisted $S^1$ over $\mathbb RP^2$.
Thus, diagrams consists of closed curves and arcs in a disk with endpoints of arcs coming in antipodal
pairs on the boundary of the disk. An example of Reidemeister moves between diagrams is presented in
Figure~\ref{fig:rp3rp3diag}.

\begin{figure}[h]
        \centering
        \subfloat{\centering
        \includegraphics[page=91]{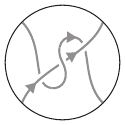}
	}        
   \raisebox{27pt}{$\;\;\xlongrightarrow{\;\Om_5\;}\;\;$}
 \subfloat{\centering
        \includegraphics[page=92]{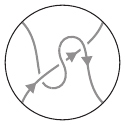}
        }
 \raisebox{27pt}{$\;\;\xlongrightarrow{\;\Om_5\;}\;\;$}
 \subfloat{\centering
        \includegraphics[page=93]{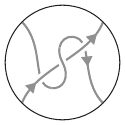}
        }
 \raisebox{27pt}{$\;\;\xlongrightarrow{\;\Om_8\;}\;\;$}
 \subfloat{\centering
        \includegraphics[page=94]{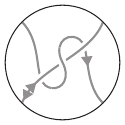}
        }
 \raisebox{27pt}{$\;\;\xlongrightarrow{\;\Om_4\;}\;\;$}
 \subfloat{\centering
        \includegraphics[page=95]{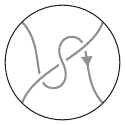}
        }
        \caption{Reidemeister moves on diagrams of a link in $\RR$}
        \label{fig:rp3rp3diag}
\end{figure}

\end{example}

\section{Arrow diagrams for links in Seifert manifolds}
\label{sec:seifert}

\begin{definition}
A standard fibered torus corresponding to a pair of coprime integers $(q,p)$, $q>0$, or $\frac{p}{q}$, is the solid cylinder $D^2 \times I $, 
where we identify the ends of the solid cylinder with a $2\pi p / q$ twist. Each $S^1$ fiber comes from $q$ vertical segments in the cylinder,
except for the core fiber which comes from the central vertical segment. We call this core fiber {\it exceptional} if $q>1$.
\end{definition}

\begin{definition}
A Seifert manifold (also a Seifert fibered space) is a closed 3-manifold which can be decomposed into a disjoint union of $S^1$'s (called fibers), 
such that each tubular neighbourhood of a fiber is a standard fibered torus.
\end{definition}


Any orientable Seifert manifold $M$ can be obtained from $F\times S^1$ or 
$N\hat{\times}S^1$ through a finite number of surgeries $(q_i,p_i)$ on 
vertical $S^1$ fibers (see \cite{O}).

\begin{figure}[h]
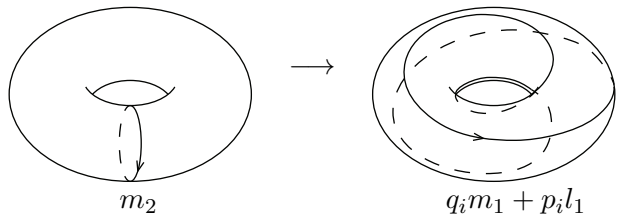

     \centering
        \begin{overpic}[page=36]{images}\put(46,-5){$m_2$}
        \end{overpic}
        \raisebox{4.2em}{$\;\;\longrightarrow\;\;$}
        \begin{overpic}[page=38]{images}\put(32,-5){$q_i m_1 + p_i l_1$}
        \end{overpic}
       \caption{Glueing map of the surgery.}\label{fig:surg}
\end{figure}

To perform such a surgery one removes a vertical solid torus $T_1$ with
longitude $l_1$ and meridian $m_1$ on $\partial T_1$ corresponding to a
vertical and horizontal $S^1$'s in the product $F\times S^1$ or 
$S_n\times S^1$ in the case of $N\hat{\times}S^1$ (see the previous section).
Then, another solid torus $T_2$ with fixed longitude $l_2$ and meridian $m_2$
is glued to $\partial T_1$, so that $m_2$ is glued to the curve $q_i m_1+p_i l_1$, see Figure~\ref{fig:surg}.
After glueing the meridional disk along this curve, 
the remaining ball of $T_1$ is glued to finish the surgery.

\begin{figure}[h]
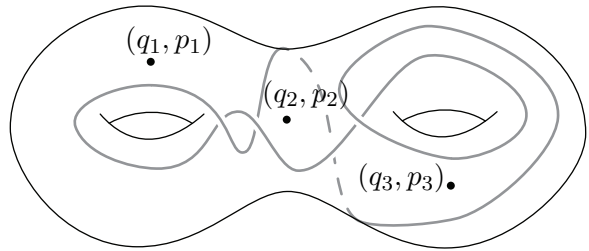

	\centering
        \begin{overpic}[page=39]{images}\put(21,35){$(q_1,p_1)$}
        \put(44,26){$(q_2,p_2)$}\put(60,12){$(q_3,p_3)$}
        \end{overpic}
       \caption{A diagram of a link in a Seifert manifold with three surgeries.}\label{fig:double}
\end{figure}

The diagram of a link $L$ in the Seifert manifold $M$ is constructed as before,
assuming that $L$ misses the exceptional fibers of the surgeries (which can be done
by general position). These exceptional fibers project to points in $F$ or $N$,
disjoint from the curves of the diagram. These points appear in the diagrams,
together with the type of surgery $(q_i,p_i)$ next to them. See Figure~\ref{fig:double} as an example.

If $q_i=1$, the fiber is not exceptional and, as a shorthand, next to the point onto which it is projected,
we put $(p_i)$ instead of $(1,p_i)$.
For instance, if $M$ is an $S^1$-bundle over $F$, then there is a unique $(p)$
fiber in the diagrams, $p\in\mathbb Z$. If $p=0$ one gets just $S_1\times F$.

\begin{figure}[h]
	\centering
		\subfloat[Forbidden position]{\hspace{1cm}\centering
	\includegraphics[page=59]{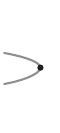}
	\label{fig:forbidslide}\hspace{1cm}
	}\hspace{0.6cm}
	\subfloat[$\Om_{(q,p)}$]{\centering
	\includegraphics[page=34]{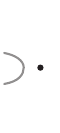}\raisebox{30pt}{$\;\;\longleftrightarrow\;\;$}\includegraphics[page=35]{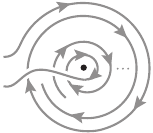}
	\label{fig:reidpq}
	}
       \caption{A forbidden position and the corresponding slide move $\Om_{(q,p)}$.}\label{fig:slide}
\end{figure}

With the added surgery fibers we get an additional forbidden position:
\begin{enumerate}[i)]\setcounter{enumi}{8}
	\item the surgery point and strand coincide (Figure~\ref{fig:forbidslide}),\hfill {(2)}
\end{enumerate}
which gives rise to the Reidemeister move $\Om_{(q_i,p_i)}$, 
corresponding to sliding an arc of the link $L$ through the meridional
disk of $T_2$. The $\Om_{(q,p)}$ move is shown in Figure~\ref{fig:slide}, it 
consists of going $q$ times around the exceptional
point and adding $p$ arrows uniformly on every $2\pi q /p$ angle.



\begin{theorem}\label{thm:reid3}
A link $L_1$ is ambient isotopic to a link $L_2$ in an orientable Seifert manifold $M$ if and only if an arrow diagram
$D_1$ of $L_1$ can be obtained from an arrow diagram $D_2$ of $L_2$ by a series of Reidemeister moves $\Om_1$ to $\Om_8$ and $\Om_{(q_i,p_i)}$.  
\end{theorem}

Let $D$ be an arrow diagram of a link $L$. We call a component of $D$ an {\it oval} if it is a component without crossing with possible arrows on it. We call an oval {\it nested} if it lies in the interior of a disk bound by another oval.

\begin{figure}[h]
	\centering
	\begin{overpic}[page=25]{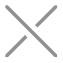}\put(37,76){\rotatebox{90}{$\bigg\}$}}\put(48,88){$n$}\end{overpic}
	\raisebox{3em}{$\;\;\sim\;\;$}
	\begin{overpic}[page=26]{skein}\put(46,82){$n$}\end{overpic}
	\raisebox{3em}{$\;\;\sim\;\;$}
	\begin{overpic}[page=27]{skein}\put(40,82){$-n$}\end{overpic}
       \caption{Oval notation.}\label{fig:on}
\end{figure}

\begin{figure}[h]
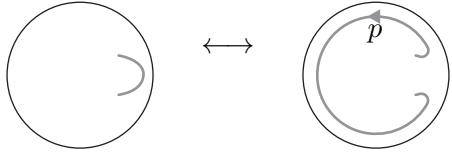

	\centering
	\begin{overpic}[page=96]{images} \end{overpic}
	    \raisebox{3.5em}{$\;\;\;\;\longleftrightarrow\;\;\;\;$}
       \begin{overpic}[page=97]{images} \put(45,74){$p$}\end{overpic}
	   \caption{The $\Omega_{(p)}$ move.}\label{fig:slidedisk}
\end{figure}

If $n$ arrows lie consecutively we simplify the diagram by placing only one arrow with an integer $n$ next to it. We interpret a negative integer above an arrow as $|n|$ reversed arrows, see Figure~\ref{fig:on}.
In lens spaces we can alternatively think of $L(p,q)$ as a manifold with Heegaard genus 1 decomposition of two solid tori. 
We isotope $L$ into the first solid torus and project it to annulus as before and an arrow diagram of a link in $L(p,q)$ can thus be viewed as a diagram on a disk. 
For $L(p,1)$, the move $\Omega_{(1,p)}$, which we will denote by $\Omega_{(p)}$ (see \cite{M1}), 
is a winding around the boundary of the disk with $p$ arrows added, see Figure~\ref{fig:slidedisk}.

\subsection{Alternative diagrams for links in Seifert manifolds}

On some occasions it may be convenient to cut the base surface to its fundamental polygon and get diagrams on a regular $n$-gon. 
Such diagrams can be expanded to all Seifert manifolds, orientable and non-orientable, 
since any Seifert manifold can be obtained from an $S^1$-bundle over $F$ or $N$ through a finite number of surgeries $(q_i, p_i)$ on vertical $S^1$ fibers.
These diagrams were introduced in \cite{M2}, see also \cite{DLP}.

\begin{figure}[h]
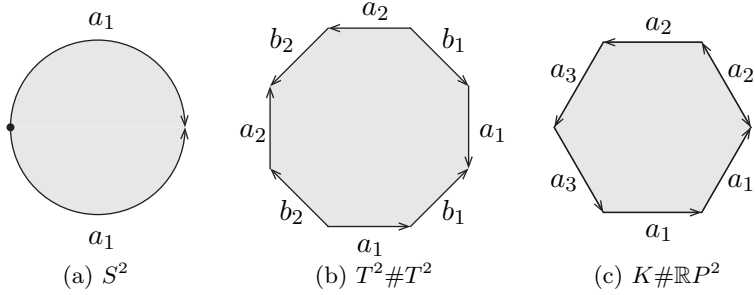

\centering
\subfloat[$S^2$]{\centering
\begin{overpic}[page=42]{images}
\put(47,3){$a_1$}\put(47,91){$a_1$}
\end{overpic}
}\hspace{.5em}
\subfloat[$T^2\#T^2$]{\centering
\begin{overpic}[page=43]{images}
\put(14,13){$b_2$}\put(-2,46){$a_2$}\put(10,82){$b_2$}\put(47,95){$a_2$}
\put(79,82){$b_1$}\put(94,46){$a_1$}\put(79,13){$b_1$}\put(45,0){$a_1$}
\end{overpic}
}\hspace{.5em}
\subfloat[$K\#\RP^2$]{\centering
\begin{overpic}[page=44]{images}
\put(9,26){$a_3$}\put(9,70){$a_3$}\put(47,89){$a_2$}\put(80,70){$a_2$}\put(80,26){$a_1$}\put(47,6){$a_1$}
\end{overpic} 
}
\caption{Fundamental polygons, where $K$ is the Klein bottle.}
\label{fig:fundamental}
\end{figure}

We start by taking the fundamental polygon $G$ of the surface, 
with the standard identification of the edges of $G$. 
We distiguish between three cases: either the surface is $S^2$, a genus $g>0$ surface $F$, or a non-orientable surface $N$, 
(see Figure~\ref{fig:fundamental}). 

Take $G \times [0,1]$. By glueing $\{x\}\times \{  0 \}$ to $\{x\}\times \{  1 \}$ for each $x\in G$, we get the trivial circle bundle $G \times S^1$. 
Since $G$ is a disk, we can orient all the fibers $\{x\} \times S^1$ coherently. If two oriented edges $a_i$ and $a_i'$ are identified in $G$, in order to get $F$, 
we can identify the cylinders $a_i \times S^1$ and $a_i' \times S^1$ in two essentially different ways: 
$a_i \times S^1$ can be glued to $a_i' \times S^1$ by identity or by a reflection on the $S^1$ component.
We assign to each edge sign $\pm1$, which takes the value $+1$ if the identification was made by identity 
and $-1$ if the identification was made by reflection. 

After the above identifications the resulting space is an $S^1$-bundle over $F$.
Any Seifert manifold can be obtained from this $S^1$-bundle by performing a finite number of $(q_i, p_i)$-surgeries on vertical fibers.




We remark that for an orientable Seifert fibered space $M$, the base space does not need to be oriented, but the edge signs are determined. 



Since the vertical projection maps, as before, an exceptional fiber to a point in the base space, it is enough to specify the image of each exceptional fiber in $G$, 
which is done by placing a point on $G$ decorated by the surgery coefficient $(q_i, p_i)$ of the fiber. 

\begin{figure}[htb]
	\centering
	\subfloat[]{\centering
	\includegraphics[page=50]{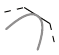}\label{fig:forbid9}
	}\hspace{0.8cm}
	\subfloat[]{\centering
	\includegraphics[page=51]{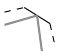}\label{fig:forbid10}
	}\hspace{0.8cm}
	\subfloat[]{\centering
	\includegraphics[page=52]{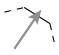}\label{fig:forbid11}
	}\hspace{0.8cm}
	\subfloat[]{\centering
	\includegraphics[page=53]{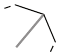}\label{fig:forbid12}
	}
	\caption{Forbidden positions of regular diagrams.}
  	\label{fig:forbidden3}
\end{figure}

We call a diagram regular if, in addition to forbidden positions (1) and (2), none of the following situations occur on the diagram:
\begin{enumerate}[i)]\setcounter{enumi}{9}
	\item border tangency~(Fig.~\ref{fig:forbid9}),
	\item crossing lies on the border~(Fig.~\ref{fig:forbid10}),
	\item arrow lies on the border~(Fig.~\ref{fig:forbid11}),
	\item arc goes through the basepoint (the preimage of the 0-cell of $F$)~(Fig.~\ref{fig:forbid12}).
\end{enumerate}

\begin{figure}[h]
	\centering
	
	\subfloat[$\Om_9$]{\centering
	\includegraphics[page=20]{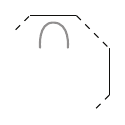}\raisebox{26pt}{$\longleftrightarrow$}\includegraphics[page=21]{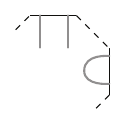}\	
	\label{fig:reid9}}\hspace{0.8cm}
	\subfloat[$\Om_{10}^{O+}$]{\centering
	\includegraphics[page=22]{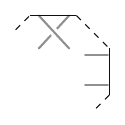}\raisebox{26pt}{$\longleftrightarrow$}\includegraphics[page=23]{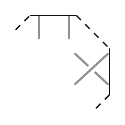}\	
	\label{fig:reid10op}}\hspace{0.8cm}
	\subfloat[$\Om_{10}^{O-}$]{\centering
	\includegraphics[page=22]{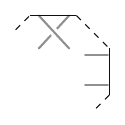}\raisebox{26pt}{$\longleftrightarrow$}\includegraphics[page=60]{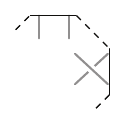}\	
	\label{fig:reid10om}}\hspace{0.8cm}
	\subfloat[$\Om_{10}^{N+}$]{\centering
	\includegraphics[page=24]{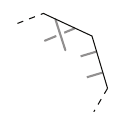}\raisebox{26pt}{$\longleftrightarrow$}\includegraphics[page=61]{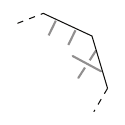}\	
	\label{fig:reid10np}}\hspace{0.8cm}
	\subfloat[$\Om_{10}^{N-}$]{\centering
	\includegraphics[page=24]{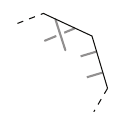}\raisebox{26pt}{$\longleftrightarrow$}\includegraphics[page=25]{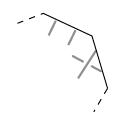}\	
	\label{fig:reid10nm}}\hspace{0.8cm}
	\subfloat[$\Om_{11}^+$]{\centering
	\includegraphics[page=30]{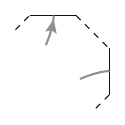}\raisebox{26pt}{$\longleftrightarrow$}\includegraphics[page=31]{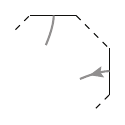}\	
	\label{fig:reid11p}}\hspace{0.8cm}
	\subfloat[$\Om_{11}^-$]{\centering
	\includegraphics[page=32]{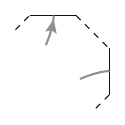}\raisebox{26pt}{$\longleftrightarrow$}\includegraphics[page=33]{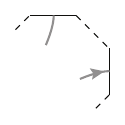}\
	\label{fig:reid11m}}\hspace{0.8cm}
	\subfloat[$\Om_{12}^O$]{\centering
	\includegraphics[page=26]{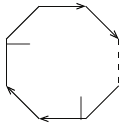}\raisebox{26pt}{$\longleftrightarrow$}\includegraphics[page=27]{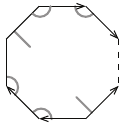}\	
	\label{fig:reid12}}\hspace{0.8cm}
	\subfloat[$\Om_{12}^S$]{\centering
	\includegraphics[page=40]{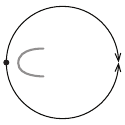}\raisebox{26pt}{$\longleftrightarrow$}\includegraphics[page=41]{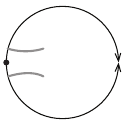}\	
	\label{fig:reid12S}}\hspace{0.8cm}
	\subfloat[$\Om_{12}^N$]{\centering
	\includegraphics[page=28]{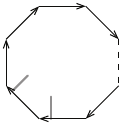}\raisebox{26pt}{$\longleftrightarrow$}\includegraphics[page=29]{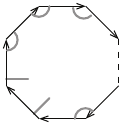}\	
	\label{fig:reid12N}}
	\caption{Additional Reidemeister moves.}
	\label{fig:reidemeister2}
\end{figure}

Positions x), xi), and xii) generate the Reidemeister moves $\Om_9$, $\Om_{10}^{O/N \pm}$, and $\Om_{11}^\pm$, that act across edges in $G$ (Figure~\ref{fig:reidemeister2}). 
The move $\Om_{10}$ comes in four flavours: the base surface is orientable (O) or non-orientable (N) and the sign of the edge is positive or negative. 
Similarly, the sign of the move $\Om_{11}^\pm$ corresponds to the sign of the edge we are pushing the arrow through. 

\begin{figure}[htb]
	\centering
	\begin{overpic}[page=37,scale=0.9]{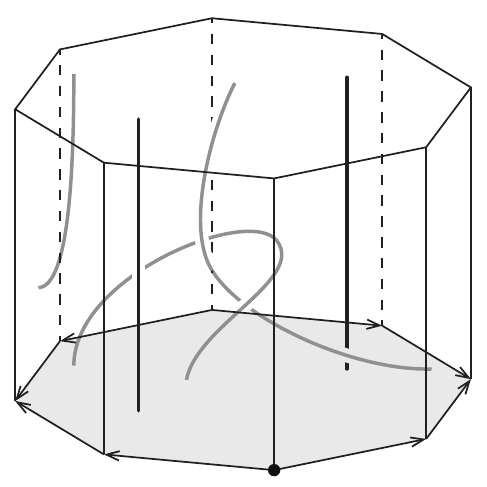}
		\put(21,35){$a_1$}\put(71,35){$a_2$}
		\put(20,6){$b_1$}\put(76,6){$b_2$}
	\end{overpic}\raisebox{34pt}{$\;\longleftrightarrow\;$}\begin{overpic}[page=39,scale=0.9]{seif}
	\put(21,35){$a_1$}\put(71,35){$a_2$}
	\put(20,6){$b_1$}\put(76,6){$b_2$}
\end{overpic}
\caption{Visualization of the move $\Om_{12}^O$.}
\label{fig:vis}
\end{figure}

Position xiii) generates the Reidemeister move $\Om_{12}$ that tells us what happens when we push an arc over the basepoint. 
The move comes in three flavours: if $G$ is a orientable genus $g>0$ surface we have $\Om_{12}^O$, if $G$ is the 2-sphere we have $\Om_{12}^S$, 
and if $G$ is a non-orientable surface we have $\Om_{12}^N$. Figure~\ref{fig:vis} shows the geometrical interpretation of $\Om_{12}$ in the case of a double torus. 

Considering the arguments above, we can now formulate the following Reidemeister theorem for links in Seifert manifolds.
\begin{theorem}
Two arrow diagrams for links in a Seifert manifold $M$ represent the same link up to ambient isotopy if and only if 
they are connected through a finite series of Reidemeister moves $\Om_1$ -- $\Om_5$, $\Om_9$ -- $\Om_{12}$,and $\Om_{(q_i,p_i)}$.
\end{theorem}




\section{Diagram conversions for links in the solid torus} 
\label{sec:passing}

In this section we define two widely used diagrams of links in the solid torus and links in lens spaces, namely the classical diagram (e.g. \cite{T,GM}) 
and the mixed link diagram (e.g. \cite{HK, HP, LR}). In addition we show how to pass between these diagrams.

\subsection{Classical diagrams}

\begin{figure}[h]
	\centering
	   \includegraphics[page=50]{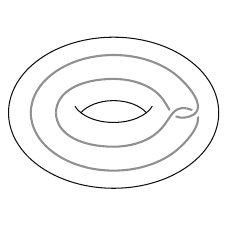}
	    \raisebox{6em}{$\;\;\;\;\longrightarrow\;\;\;\;$}
       \includegraphics[page=51]{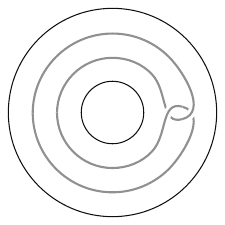}
	   \caption{A link in the solid torus (left) and its classical diagram (right).}\label{fig:torus}
\end{figure}

Let $T = D^2 \times S^1$ be the solid torus and let $L \subset T$ be a link in $T$. 
A classical link diagram of $L$ is the projection of $L$ to the annulus in $T$ that lies in the plane spanned by the longitude $l$ of $T$, see Figure~\ref{fig:torus}.

As in the previous section, we can think of $L(p,q)$ as glueing of two solid tori.  
A classical diagram of a link in $L(p,q)$ can again be viewed as a diagram of a link in a solid torus.  
Alternatively, one can think of $L(p,q)$ as the result of a $p/q$ rational surgery performed on a unknot $C$ in $S^3$, 
where again we can project the knot to the annulus $S^3 \setminus C$


\subsection{Passing from classical to arrow diagrams}
\begin{figure}[h]
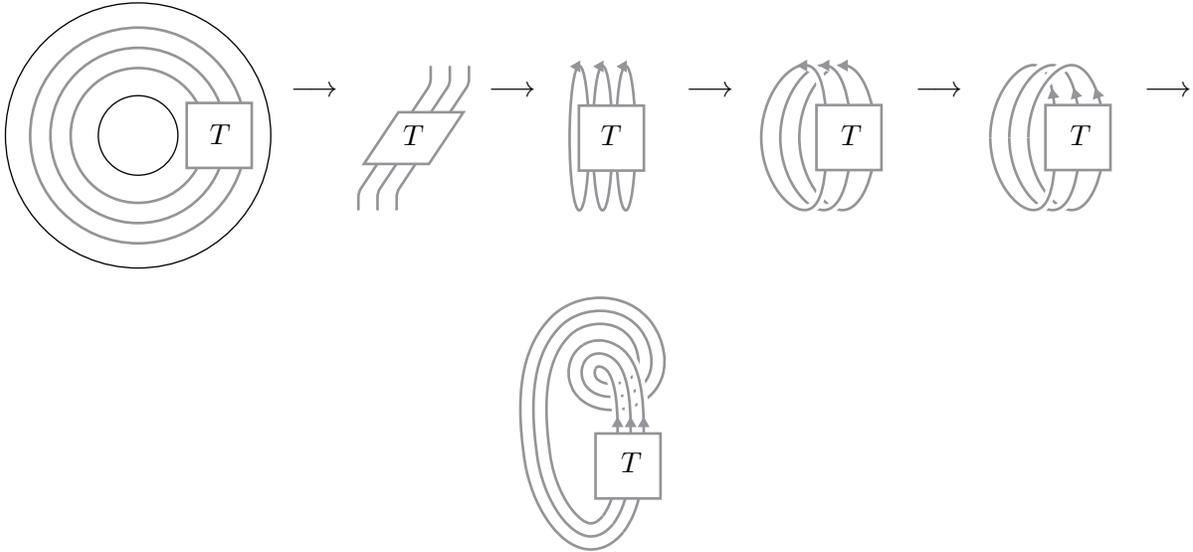

	\centering
	\begin{overpic}[page=53]{images}\put(75,47){$T$}\end{overpic}
	   \raisebox{6.3em}{$\longrightarrow$}
      \raisebox{1.0em}{\begin{overpic}[page=54]{images}\put(25,47){$T$}\end{overpic}}
       \raisebox{6.3em}{$\longrightarrow$}
      \raisebox{1.0em}{\begin{overpic}[page=55]{images}\put(25,47){$T$}\end{overpic}}
       \raisebox{6.3em}{$\longrightarrow$}
      \raisebox{1.0em}{\begin{overpic}[page=56]{images}\put(44,47){$T$}\end{overpic}}
       \raisebox{6.3em}{$\longrightarrow$}
      \raisebox{1.0em}{\begin{overpic}[page=57]{images}\put(44,47){$T$}\end{overpic}}
       \raisebox{6.3em}{$\longrightarrow$}
      \raisebox{1.0em}{\begin{overpic}[page=58]{images}\put(41,32){$T$}\end{overpic}}
       \caption{Passing from classical diagrams to arrow diagrams.}\label{fig:class_to_arrow}
\end{figure}

Consider a classical diagram of a link in $D\times S^1$. Each 
such diagram  can be obtained by closing a $(n,n)$-tangle $T$ with $n$ 
strands parallel to $\{P\}\times S^1$ for any point $P \in D$. 
The construction of the corresponding arrow diagram is presented in 
Figure~\ref{fig:class_to_arrow}. By rotating the tangle 
it can be made horizontal, so that in the arrow diagram it will look 
also like $T$. In the arrow diagram, the strands become arcs starting 
at the upper endpoints of the tangle, with arrows on them, then going 
under the tangle and joining the lower endpoints of the tangle.
Applying some $\Om_5$ moves all arrows can be moved to the upper 
endpoints of the tangle. Then one notices that the strands go from the 
upper to the lower endpoints of $T$ with a full negative twist.

Summing up, if one has a classical diagram of a link in the solid torus, 
presented as a closure of a tangle, an arrow diagram for this link 
is obtained from the same tangle, by adding arrows going up at the upper 
endpoints of the tangle, making a full negative twist of the strands and 
closing the tangle {\it on the left}. By a similar construction, one can 
add a positive twist and close the tangle on the right.




\subsection{Passing from arrow to classical diagrams}
\begin{figure}[h]
	\centering
	   \includegraphics[page=59]{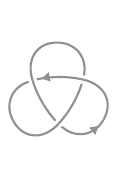}
	   \raisebox{5em}{$\;\;\longrightarrow\;\;$}
       \includegraphics[page=60]{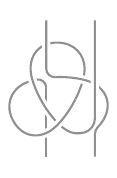}
       \raisebox{5em}{$\;\;\longrightarrow\;\;$}
       \includegraphics[page=61]{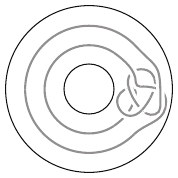}
       \caption{Passing from arrow diagrams to classical diagrams.}\label{fig:arrow_to_class}
\end{figure}

An arrow diagram can be viewed as an almost flat diagram outside 
small neighborhoods of the arrows (i.e. it lies in a thickened 
$D\times\{1\}$ in $D\times S^1$). The neighborhoods of arrows correspond 
to vertical 
strands parallel to $\{P\}\times S^1$, $P\in D$. We choose an arbitrary 
direction in $\mathbb R^2$ of the arrow diagram (for instance the 
vertical one), and rotate the diagram around the axis orthogonal to this
direction.
One may assume, by general position, that the arrows do not point in the 
chosen direction. Then the arrows become vertical strands: just before 
the arrow a vertical strand goes up above other strands, and just after 
the arrow a vertical strand goes from below under other strands. Closing 
these vertical strands in an annulus gives a classical diagram from the 
original arrow diagram. An example is presented in 
Figure~\ref{fig:arrow_to_class}.

\subsection{Mixed link diagrams}

Let $N(U)$ be a the thickened unknot $U$ in $S^3$. Since $T = S^3 \setminus N(U)$ is a solid torus, we can represent any link $L$ in $T$ with a diagram of $L \cup U$ in the plane. We call $L$ the moving component and $U$ the fixed component. We also keep track of the two types of component, by coloring them with two different colors. Such a diagram is called a mixed link diagram for $T$.

\begin{figure}[h]
	\centering
	   \includegraphics[page=62]{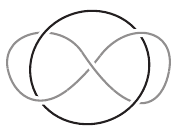}
       \caption{A mixed link diagram for the solid torus.}\label{fig:mixed1}
\end{figure}

Every closed oriented 3-manifold $M$ can be constructed from a link $A \subset S^3$ 
on which we perform (integer or rational) Dehn surgeries on its components~\cite{Lick,Wal}. 
Furthermore, each component of the link $A$ can be assumed to be unknotted. In this way, we can represent a link $L$ in $M$ by a diagram $A \cup L$ in the plane, 
but again, we keep track of the fixed and moving parts by coloring them with two distinct colors; 
in addition, we equip each component of $A$ with the surgery coefficient.

\begin{figure}[h]
	\centering
	\begin{overpic}[page=62]{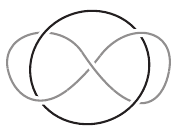}
		\put(60,-3){$(q,p)$}\end{overpic}
      \caption{A mixed link diagram for $L(p,q)$.}\label{fig:mixed2}
\end{figure}

For example, the lens space $L_{p,q}$ is the result of a $(q,p)$ or $\frac{p}{q}$ rational surgery on the unknot in $S^3$.
Figure~\ref{fig:mixed2} shows an example of a knot in $L(p,q)$ (see \cite{DL1, LR1, LR2}).

\subsection{Passing from classical diagrams to mixed link diagrams for the solid torus}

\begin{figure}[h]
	\centering
	   \includegraphics[page=51]{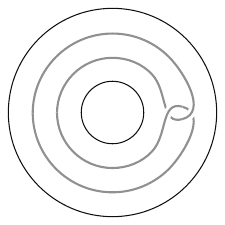}
	   \raisebox{6em}{$\;\;\longrightarrow\;\;$}
       \includegraphics[page=63]{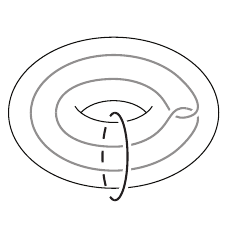}
       \raisebox{6em}{$\;\;\longrightarrow\;\;$}
       \includegraphics[page=100]{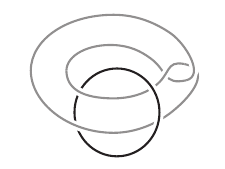}
       \caption{Passing from classical diagrams to mixed link diagrams.}\label{fig:class_to_mixed}
\end{figure}

Passing from classical diagrams to mixed link diagrams and back is easy: 
the complement of the solid torus is the thickened $N(U)$, thus $U$ represents the fixed component, see Figure~\ref{fig:class_to_mixed}.

\subsection{Passing from mixed link diagrams to classical diagrams for the solid torus}

\begin{figure}[h]
	\centering
	   \includegraphics[page=65]{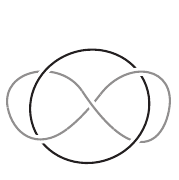}
	   \raisebox{4em}{$\;\;\longrightarrow\;\;$}
       \includegraphics[page=68]{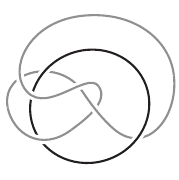}
       \raisebox{4em}{$\;\;\longrightarrow\;\;$}
       \includegraphics[page=69]{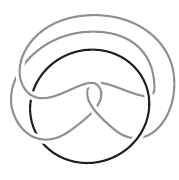}
       \raisebox{5em}{$\;\;\longrightarrow\;\;$}
       \includegraphics[page=70]{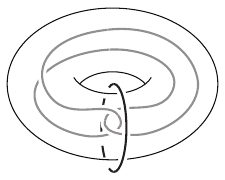}
      \raisebox{5em}{$\;\;\longrightarrow\;\;$}
       \includegraphics[page=71]{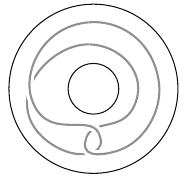}
       \caption{Passing from mixed link diagrams to classical diagrams.}\label{fig:mixed_to_class}
\end{figure}

To pass from mixed link diagrams to classical diagrams, we isotope the moving components of mixed link in such a way, that overcrossings with the fixed component $U$ lie on one side, say on the left and undercrossings lie on the other side. 
Strands connecting undercrossings and overcrossings should also connect on one side, say the top, see Figure~\ref{fig:mixed_to_class}.

\section{The Kauffman bracket and HOMFLYPT skein modules}
\label{sec:skein_modules}

Let $M$ be an orientable 3-manifold.
Take a coefficient ring $R$ and a unit $A \in R$ (an element with a multiplicative inverse). And let $\calL_{\mathrm{fr}}(M)$ be the set of isotopy classes of framed links in $M$, including the class of the empty link $[\emptyset]$. Let $R\calL_{\mathrm{fr}}(M)$ be the free $R$-module spanned by $\calL_{\mathrm{fr}}(M)$.

We would like to impose the Kauffman relation and the framing relation in $R\calL_{\mathrm{fr}}(M)$. We therefore take the submodule $\mathcal{S}_{\mathrm{fr}}(M)$ of $R\calL_{\mathrm{fr}}(M)$ generated by
\begin{align*}
\tag{Kauffman relator} \raisebox{-4.5pt}{\includegraphics[page=1]{skein}} - A \raisebox{-4.5pt}{\includegraphics[page=2]{skein}} -A^{-1} \raisebox{-4.5pt}{\includegraphics[page=3]{skein}},\\
\tag{framing relator} L \sqcup \raisebox{-4.5pt}{\includegraphics[page=4]{skein}} - (-A^2 - A^{-2}) L.
\end{align*}

The Kauffman bracket skein module $\kbsm(M)$ is defined as 
$R\calL_{\mathrm{fr}}(M)$ modulo these two relations: $$\kbsm(M) = R\calL_{\mathrm{fr}}(M) / \calS(M).$$ 

\begin{example}
For the 3-sphere, $\kbsm(S^3)$ is a free $R$-module with a basis consiting of a single element, the equivalence class of the unknot (here, we exclude the empty link). 
Expressing a link in this basis and, for the unknot, evaluating $[O] = 1$, we get exactly the Kauffman bracket.
\end{example} 







\begin{figure}[h]
	\centering
	\includegraphics[page=75]{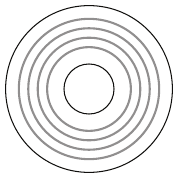}\raisebox{4.8em}{$\;\;\sim\;\;$}\includegraphics[page=76]{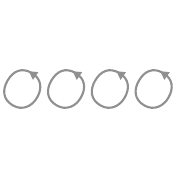}
\caption{Generator $x^4$}\label{fig:xn}
\end{figure}

\begin{theorem}[\cite{T}]\label{thm:kbsmt}
Let $T$ be the solid torus, $\kbsm(T)$ is freely generated by an infinite set of generators $\{x^n\}_{n=0}^{\infty}$, 
where $x^n$, $n>0$, is a parallel copy of $n$ longitudes of $T$ (see Figure~\ref{fig:xn}) and $x^0$ is the empty link.
\end{theorem}

\begin{theorem}[\cite{HP}]\label{thm:kbsmlpq}
$\kbsm(L(p,q))$ is freely generated by the set of generators $\{x^n\}_{n=0}^{\lfloor p/2 \rfloor}$, 
where $x^n$, $n>0$, is a parallel copy of $n$ longitudes of $T \subset L(p,q)$ (see Figure~\ref{fig:xn}) and $x^0$ is the empty link.
\end{theorem}

For the HOMFLYPT skein module we take oriented (unframed) links and impose on them the HOMFLYPT skein relation.

Let the ring $R$ have two units $v, z \in R$. Let $\calL_{\mathrm{or}}(M)$ be the set of isotopy classes of oriented links in $M$, including the class of the empty link $[\emptyset]$ and let $R\calL_{\mathrm{or}}(M)$ be the free $R$-module spanned by $\calL_{\mathrm{or}}(M)$. 

We take the submodule $\calS(M)$ of $R\calL_{\mathrm{or}}(M)$ generated by the expressions 
\begin{align*}
\tag{HOMFLYPT relator} v^{-1} \raisebox{-4.5pt}{\includegraphics[page=21]{skein}} - v \raisebox{-4.5pt}{\includegraphics[page=22]{skein}} - z \raisebox{-4.5pt}{\includegraphics[page=23]{skein}}.
\end{align*}

We also add to $\calS(M)$ the HOMFLYPT relation involving the empty knot:\begin{align*}
\tag{empty knot relator} v^{-1}\emptyset - v\emptyset - z  \raisebox{-4.5pt}{\includegraphics[page=4]{skein}}.
\end{align*}

The HOMFLYPT skein module $\hsm(M)$ of $M$ is $R\calL_{\mathrm{or}}(M)$ modulo the above relations:
$$\hsm(M) = R\calL(M)/\calS(M).$$


\begin{figure}[h]
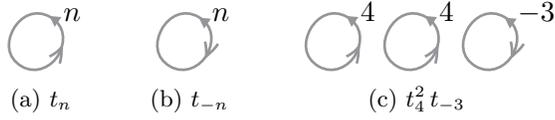

	\centering
	\subfloat[$t_n$]{\centering
		\begin{overpic}[page=86]{images}\put(90,80){$n$}\end{overpic}
		\label{fig:tn}
	}\hspace{0.8cm}
	\subfloat[$t_{-n}$]{\centering
		\begin{overpic}[page=88]{images}\put(90,80){$n$}\end{overpic}
		\label{fig:tmn}
	}\hspace{0.8cm}
	\subfloat[$t_4^2\,t_{-3}$]{\centering
		\begin{overpic}[page=86]{images}\put(90,80){$4$}\end{overpic}
		\begin{overpic}[page=86]{images}\put(90,80){$4$}\end{overpic}
		\begin{overpic}[page=88]{images}\put(90,80){$-3$}\end{overpic}
		\label{fig:ttt}
	}
	\caption{HOMFLYPT skein module generators.}
	\label{fig:t}
\end{figure}


Let $t_n$, $n \in \Z$ be the knot in Figures \ref{fig:tn} and \ref{fig:tmn}, note that for $n < 0$, $t_n$ is $t_{|n|}$ with reversed orientation.

\begin{theorem}[\cite{HK, T}]
The HOMFLYPT skein module of the solid torus $T$ is a free $R$-module, generated by the infinite set $$\calB = \{ t_{k_1} \cdots t_{k_s} \;|\;k_i \in \Z \setminus \{0\},\; k_1 \leq \cdots \leq k_s  \} \cup \{\emptyset\}.$$
\end{theorem}

\subsection{Results for Kauffman bracket and HOMFLYPT skein modules using arrow diagrams}

The first result, using arrow diagrams, was to compute the Kauffman bracket skein module of a product
of $S^1$ with a disk with two holes. There is an incompressible, non-boundary parallel torus immersed
(not embedded) in such a manifold. In several examples of manifolds, such surfaces, when embedded, yield torsion,
so it was interesting to consider a case with an immersion instead of an embedding. There was no torsion in
this case:

\begin{theorem}[\cite{MD}]\label{thm:disk_two_holes}
Let $M$ be the product of a disk with two holes and $S^1$. Then $\kbsm(M)$ is freely generated by an infinite set of generators.
\end{theorem}

Now, consider links in the manifold $\RR$ (see Example~\ref{ex:rp3rp3}).
Let $t=-A^{-3}x$ (see Figure~\ref{fig:xn}).
We can view $t$ as an oval with one arrow, obtained from $x$ by adding a negative kink.
The multiplication of two arrow diagrams of links in the solid torus consists in putting them in two disjoint disks.
Thus, for example, $x^3=x^2 x=xxx$ is the diagram were the three $x$'s are in three disjoint disks.

Let $Q_n$, $n\in\N\cup\{0\}$ be defined by:
\begin{equation*}Q_0=1,\; Q_{1}=t\; and\;Q_n=t Q_{n-1}-Q_{n-2}
\end{equation*}

\begin{figure}[h]
\centering
	\subfloat[$E$]{\centering
	\includegraphics[page=85]{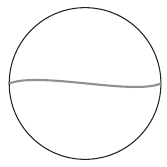}
	\label{fig:e}
	}\hspace{0.8cm}
	\subfloat[$E'$]{\centering
	\includegraphics[page=84]{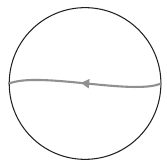}
	\label{fig:ep}
	}
	\caption{Generators $E$ and $E'$.}
	\label{fig:es}
\end{figure}

Let $E$ be the knot with diagram in Figure~\ref{fig:e}: an arc with two antipodal endpoints and no crossings.
Let $E'$ be the knot with diagram in Figure~\ref{fig:ep}: an arc as in $E$ with an added arrow.
\begin{theorem}[\cite{M1}]\label{thm:rp3rp3}
$\KB(\RR)=R\oplus R\oplus R[t]/S$, where $R=\Z[A,A^{-1}]$ and $S$ is the submodule
of $R[t]$ generated by:

$(A^{n+1}+A^{n-1})(Q_n-1)-2(A+A^{-1})\Sum_{k=1}^{\frac{n}{2}}
A^{n+2-4k}$, for $n\ge 2$ even,

$(A^{n+1}+A^{n-1})(Q_n-t)-2t\Sum_{k=1}^{\frac{n-1}{2}} A^{n+1-4k}$, for
$n\ge 3$ odd.

The generators of the two $R$'s are $E$ and $E'$.
\end{theorem}

By adding $(q,p)$ fibers to $\RR$ one gets the prism manifolds. If $q=1$, the prism manifold 
is denoted by $M_p$ (it has no exceptional fibers, just as $L(p,1)$).

\begin{theorem}[\cite{M2}]\label{thm:prisms}
$\KB(M_p)$ is a free $R$-module generated by $\emptyset$, $x$,
$x^2$,\ldots,$x^{1+\lfloor\frac{p}{2}\rfloor}$, $E$ and (if $p$ is
even) $E'$.
Thus, it has $3+\lfloor\frac{p}{2}\rfloor$ generators, if $p$ is odd, 
and $4+\frac{p}{2}$ generators, if $p$ is even.
\end{theorem}




Let us now turn to HOMFLYPT skein modules. 
Consider now links in the lens space $L_{p,1}$. Recall, from Section 2, that  
arrow diagrams of such links lie in the disk and there is an additional slide move $\Om_{(p)}$ (Figure \ref{fig:slidedisk}). 

Let $$\mathcal{B}_p=\{t_{k_1} \cdots t_{k_s} \; | \; k_i\in\Z \setminus \{0\},\;
-\frac{p}{2} < k_1 \leq \cdots \leq k_s \leq \frac{p}{2}\}\cup\{\emptyset\}.$$
Using the move $\Om_{(p)}$, elements in $\mathcal{B}$ can be expressed with elements
in $\mathcal{B}_p$. In fact, one gets more:

\begin{theorem}[\cite{GM}]\label{thm:homflypt_lp1}
$\hsm(L_{p,1})$ is free with basis $\mathcal{B}_p$.
\end{theorem}

See also \cite{DLP} for the braid approach to the HOMFLYPT skein module of $L(p, 1)$.





\section{Alternative bases for the Kauffman bracket skein modules}
\label{sec:alternative_kbsm}

\begin{figure}\centering
\subfloat[$P_4$]{\centering
	\includegraphics[page=77]{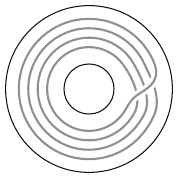}\raisebox{4.8em}{$\;\;\sim\;\;$}\begin{overpic}[page=78]{images}\put(22,61){$4$}\end{overpic}
	\label{fig:std2}
	}\hspace{1.2cm}
	\subfloat[$P_{-4}$]{\centering
	\includegraphics[page=79]{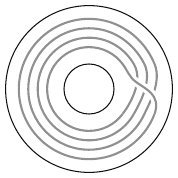}\raisebox{4.8em}{$\;\;\sim\;\;$}\begin{overpic}[page=80]{images}\put(22,61){$4$}\end{overpic}
	\label{fig:std3}
	}\hspace{0.8cm}
	\subfloat[$y_4$]{\centering
	\includegraphics[page=82]{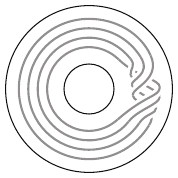}\raisebox{4.8em}{$\;\;\sim\;\;$}\begin{overpic}[page=83]{images}\end{overpic}
	\label{fig:std4}
	}
	\caption{Diagrams of $P_4$, $P_{-4}$, and $y_4$}
	\label{fig:std}
\end{figure}

We denote by $P_n$ an oval with $n$ counterclockwise arrows ($|n|$ clockwise arrows if $n<0$) 
and by $y_n$ a nested system of $n$ ovals with one counterclockwise arrow on each if $n>0$ or clockwise arrow if $n<0$. By convention $y_0$ is the empty link.
See Figure~\ref{fig:std}.

We exhibit some alternative bases for the Kauffman bracket skein modules of the solid torus and lens spaces. 
The results follow easily from the following lemmas.

\begin{lemma}[\cite{MD}] \label{lemma:turn_one_oval} In $\kbsm(T)$ one can revert $x$ in the sense that
$$x = P_1 = A^6 \, P_{-1}.$$
\end{lemma}




\begin{lemma}[\cite{MD}]\label{lemma:turn_oval}
From an oval $P_n$, $n > 0$, we can push off an arrow using Reidemeister moves and skein relations in the sense that
$$P_n = -A^{-2} P_{n-1}x - A^2 P_{n-2}.$$ Similarly for $P_n$, $n<0$, it holds $$P_n = -A^{2} P_{n+1} P_{-1} - A^{-2} P_{n+2}.$$
\end{lemma}

The following lemma illustrates the methods used when doing calculations in $\kbsm(T)$.

\begin{lemma}\label{lemma:push_oval}
An oval can be pushed to an adjacent arc in $\kbsm(T)$ using skein relations in the sense that
$$\begin{overpic}[page=46]{skein}\put(33,90){$n$}\end{overpic}\
\raisebox{1.0em}{$=\sum_{i=-n}^{n} r_i$}\
\begin{overpic}[page=47]{skein}\put(40,40){$i$}\end{overpic}$$
for $n>0$ and $r_i \in R$, furthermore it holds that $r_n = -A^{2n+2}$ and $r_{-n} = -A^2$.

\begin{proof}

First we make an $\Om_2$ move through the arc and push an arrow through, then we resolve the newly changed crossing by the Kauffman skein relation.
For the first term we resolve the remaining crossing and for the second term we push the arrow through the remaining crossing and perform an $\Om_1$:
 
$$\begin{overpic}[page=37]{skein}\put(16,78){$n$}\end{overpic}\
\raisebox{2.4em}{$\xlongequal{\Om_2}$}\begin{overpic}[page=38]{skein}\put(16,78){$n$}\end{overpic}\
\raisebox{2.4em}{$\xlongequal{\Om_5}$}\begin{overpic}[page=39]{skein}\put(0,78){$n\!-\!1$}\end{overpic}\
\
\raisebox{2.4em}{$\xlongequal{}A$}\begin{overpic}[page=43]{skein}\put(0,78){$n\!-\!1$}\end{overpic}\
\raisebox{2.4em}{$+A^{-1}$}\begin{overpic}[page=40]{skein}\put(0,78){$n\!-\!1$}\end{overpic}$$
$$
\raisebox{2.4em}{$\xlongequal{}A^2$}\begin{overpic}[page=44]{skein}\put(0,78){$n\!-\!1$}\end{overpic}\
\raisebox{2.4em}{$+$}\begin{overpic}[page=45]{skein}\put(0,78){$n\!-\!1$}\end{overpic}\
\raisebox{2.4em}{$+A^{-1}$}\begin{overpic}[page=48]{skein}\put(0,78){$n\!-\!1$}\end{overpic}$$
$$
\raisebox{2.4em}{$\xlongequal{\Om_1}A^2$}\begin{overpic}[page=50]{skein}\put(0,78){$n\!-\!1$}\end{overpic}\
\raisebox{2.4em}{$+$}\begin{overpic}[page=52]{skein}\put(0,70){$n\!-\!2$}\end{overpic}\
\raisebox{2.4em}{$-A^{2}$}\begin{overpic}[page=51]{skein}\put(27,52){$n$}\end{overpic}$$

If $n=1$, the first two lines of the equations above show that $r_1=-A^4$ and $r_{-1}=-A^2$.
For $n>1$, we repeat this process with the oval with $n-1$ arrows to the left of the strand with one arrow up.
Then, we iterate this process with ovals with less and less arrows until there is only one arrow left. 
At each step an arrow will be transferred from the oval to the strand upwards with coefficient $A^2$.
The coefficient of $r_n$ is thus $(A^2)^{n-1}(-A^4)=-A^{2n+2}$. 
It follows from the equations above that in the procedure we get $n$ arrows pointing downwards only at the first step, so $r_{-n}$ is $-A^2$.
\end{proof}
\end{lemma}



In addition to Turaev's basis of $\kbsm(T)$ in~\cite{T} and Hoste/Przytycki's basis of $\kbsm(L(p,q))$ in~\cite{HP}, 
we show in the next propositions some alternative choices for bases of these skein modules. 

\begin{proposition}\label{prop:basis1}
The set $\{P_n\}_{n=1}^{\infty}\cup\{\emptyset\}$ forms a free basis of $\kbsm(T)$.
\begin{proof}
From Lemma~\ref{lemma:turn_oval} it follows that, for $n>0$, $P_n$ is a polynomial of degree $n$ with leading
invertible coefficient $(-1)^{n+1}A^{-2n+2}$.
Thus, the $P_n$'s can be expressed with the $x^n$'s with an upper triangular matrix with invertible coefficients on the diagonal.
It follows that  $\{P_n\}_{n=1}^{\infty}\cup\{\emptyset\}$ is a basis of $\kbsm(T)$.
\end{proof}
\end{proposition}

\begin{proposition}\label{prop:basis2}
The set $\{P_{-n}\}_{n=1}^{\infty}\cup\{\emptyset\}$ forms a free basis of $\kbsm(T)$.
\begin{proof}
From Lemma~\ref{lemma:turn_one_oval}, $P_{-1}=A^{-6}P_1=A^{-6}x$.
From Lemma~\ref{lemma:turn_oval} it follows that, for $n<0$, $P_n$ is a polynomial of degree $|n|$ with leading
invertible coefficient $(-1)^{n+1}A^{-2n-8}$. The rest of the proof is the same as in the preceding proposition.
\end{proof}
\end{proposition}

\begin{proposition}\label{prop:basis3}
The set $\{y_n\}_{n=0}^{\infty}$ forms a free basis of $\kbsm(T)$.
\begin{proof}
We will show that $y_n$ is a polynomial of degree $n$ in $x$ with an invertible leading coefficient.
Then we will be done, just as in the proofs of the preceding two propositions.

The proof is by induction on $n$. Obviously it holds for $n=1$ because $y_1=x$.
It will be useful to have this more general induction hypothesis, for $k,l\ge 0$, $k+l\le n$:

$(H_{k,l}):\;x^k$ nested inside $y_l$ is a polynomial of degree $k+l$ with an invertible leading coefficient.

For instance $x$ nested inside $y_l$ is just $y_{l+1}$. We perform an inner induction on $l$. 
The hypothesis is true for $l=0$ because $x^k$ nested in $y_0$ is just $x^k$.
Suppose $(H_{k',l'})$ is true for all $k'+l'<n$ and that it is also true for all 
$x^{k'}$ nested in $y_{l'}$ where $l'<l$, $k'+l'=n$.

Now consider $x^k$ nested in $y_l$, $k+l=n$.
Use Lemma~\ref{lemma:push_oval} to push the $k$ $x$'s into the most nested oval of $y_l$.
Pushing one such $x$, the only nonzero coefficients will be $r_1=-A^4$ and $r_{-1}=-A^2$ corresponding to adding a
counterclockwise or clockwise arrow respectively to this most nested oval of $y_l$.
As originally there is one countercklockwise arrow on this oval, this oval will become $P_m$ when the $x$'s are
pushed into it, with $-k+1\le m\le k+1$.
We know that $P_m$ is a polynomial in $x$ of degree $|m|$ with an invertible leading coefficient. 
By induction on $n$, the terms with $m<k+1$ will have degree less than $n$ (some arrows are cancelled).
The remaining term, $m=k+1$, corresponds to all $x$'s being pushed as counterclockwise arrows,
yielding $(-A^4)^k$ times $P_{k+1}$ nested inside $y_{l-1}$.
Now $P_{k+1}$ is of degree $k+1$ with an invertible leading coefficient. The terms of $P_{k+1}$ of degree
less than $k+1$ nested in $y_{l-1}$ will have degree less than $n$ by induction.
Finally, the only term remaining is an invertible coefficient times $x^{k+1}$ nested in $y_{l-1}$,
which is of degree $n$ times an invertible coefficient by induction on $l$.
Thus $(H_{k,l})$ is true.
\end{proof}
\end{proposition}

\begin{proposition}\label{prop:basis4}
The set $\{y_{-n}\}_{n=0}^{\infty}$ forms a free basis of $\kbsm(T)$.
\begin{proof}
The proof mirrors that of the previous proposition, using $P_{-1}=A^{-6}x$ instead of $x$.
\end{proof}
\end{proposition}

As, for $n>0$, $P_n$, $P_{-n}$, $y_n$ and $y_{-n}$ are all polynomials of degree $n$ in $x$ with 
an invertible leading coefficient, the following proposition follows from
Theorem~\ref{thm:kbsmlpq}.

\begin{proposition}\label{prop:basis5}
The sets $\{P_n\}_{n=1}^{\lfloor p/2 \rfloor}\cup\{\emptyset\}$, $\{P_{-n}\}_{n=1}^{\lfloor p/2 \rfloor}\cup\{\emptyset\}$, $\{y_n\}_{n=0}^{\lfloor p/2 \rfloor}$, 
and $\{y_{-n}\}_{n=0}^{\lfloor p/2 \rfloor}$ are all free bases of $\kbsm(L(p,q))$.
\end{proposition}




\section{Alternative bases for the HOMFLYPT skein modules}
\label{sec:alternative_hsm}

The following two lemmas are from \cite{GM}. Using them, we will exhibit new bases for $\hsm(L_{p,1})$.
Recall that $t_n$, $n\in\mathbb Z \setminus \{0\}$, stands for an oval with $|n|$ countercklockwise arrows on it, which is oriented in
a counterclockwise way if $n>0$ and in a clockwise way otherwise.

Denote by $\bar{t}_n$, $n\in\mathbb Z \setminus \{0\}$, the oval obtained from $t_n$ by reversing all arrows and the orientation.

\begin{figure}[h]
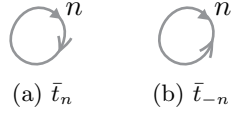

 \centering
	\subfloat[$\bar{t}_n$]{\centering
		\begin{overpic}[page=89]{images}\put(90,80){$n$}\end{overpic}
		\label{fig:btn}
	}\hspace{0.8cm}
	\subfloat[$\bar{t}_{-n}$]{\centering
		\begin{overpic}[page=87]{images}\put(90,80){$n$}\end{overpic}
		\label{fig:btmn}
	}
	\caption{Knots $\bar{t}_n$ and $\bar{t}_{-n}$ for $n>0$.}
	\label{fig:bt}
\end{figure}

\begin{lemma}\label{lemma:hsm_revert}
In $\hsm(T)$ we can revert clockwise arrows on an oval in the sense that for $n > 0$ 

$$\bar{t}_n=\Sum_i A_i T_i \mbox{  and  }\bar{t}_{-n}=\Sum_i A'_i T'_i$$

where $A_i,A'_i\in R$, $T_i,T'_i\in\mathcal{B}$.
More precisely, any $T_i=t_{k_1} \cdots t_{k_s}$, where all $k_i>0$ and $k_1+ \cdots+k_s=n$.
Similarly, any $T'_i=t_{k_1} \cdots t_{k_s}$, where all $k_i<0$ and $k_1+ \cdots+ k_s=-n$.
\end{lemma}

The following lemma is a reformulation of Lemma 2 in \cite{GM}, emphasizing the orientations and making the coefficients $A_0$ and $A'_0$ explicit. 

\begin{lemma}\label{lemma:hsm_push}

Let $D$ be a diagram of a link $L$ with an oval containing $n$ arrows, $n \in \Z$, and a strand adjacent to it, 
that may contain arrows outside the drawn region. We say that the orientations of the oval and the strand {\it agree}
if the oval has a counterclokwise orientation and the strand right to it is oriented upwards or if both have opposite orientations.
Otherwise we say that their orientations {\it disagree}.

The oval can be pushed through the strand. There are four possible configurations of orientations, but we will explicitly point out two, with $n>0$ 
and the orientations agreeing or disagreeing. The following formulas hold in $\hsm(T)$:

$$ \raisebox{-1em}{\begin{overpic}[page=54]{skein}\put(30,90){$_n$}\end{overpic}} = \sum_{i=0}^{n} A_i\;\; \raisebox{-1em}{\begin{overpic}[page=55]{skein}\put(-10,62){$_{i}$}\put(44,90){$_{n\text{-}i}$}\end{overpic}} \text{\;\;\;\;\;or\;\;\;\;\;}  \raisebox{-1em}{\begin{overpic}[page=56]{skein}\put(30,90){$_n$}\end{overpic}} = \sum_{i=0}^{n} A'_i\;\; \raisebox{-1em}{\begin{overpic}[page=57]{skein}\put(-10,62){$_{i}$}\put(44,90){$_{n\text{-}i}$}\end{overpic}},$$
where $A_i,A'_i \in R$. Furhtermore, we can keep track of the coefficients $A_0$ and $A'_0$ getting $A_0 = v^{2n}$ and $A'_0 = v^{-2n}$.
For the remaining two configurations the tranfer of arrows is as in the first formula when the orientations agree and as in the second formula
when the orientations disagree.
\end{lemma}

\begin{figure}[h]
\centering
	\subfloat[$t_{5,2,-1,-3} \in \calB'$]{\centering
		\begin{overpic}[page=58]{skein}\put(45,89){$_3$}\put(45,78.5){$_1$}\put(45,68){$_2$}\put(45,57.5){$_5$}\end{overpic}
		\label{fig:bp}
	}\hspace{0.8cm}
	\subfloat[$t_{-1,-3,2,4} \in \calB''$]{\centering
		\begin{overpic}[page=59]{skein}\put(45,89){$_4$}\put(45,78.5){$_2$}\put(45,68){$_3$}\put(45,57.5){$_1$}\end{overpic}
		\label{fig:bpp}
	}
	\caption{An element from $\calB'$ and $\calB''$.}
	\label{fig:b}
\end{figure}

In \cite{DL} a new basis for $\hsm(T)$ is presented. Translating this basis into arrow diagrams
gives a basis $\calB'$. Its elements have diagrams consisting of $s\in\mathbb N$ concentric ovals 
with $k_i\in\mathbb Z \setminus \{0\}$ arrows on each oval, denoted $t_{k_1,k_2, \ldots ,k_s}$, satisfying $k_1\ge k_2\ge \ldots \ge k_s$,
with $k_1$ arrows on the most nested oval, $k_2$ arrows on the next one and, so on, with $k_s$ arrows on the oval containing 
all other ovals, see Figure~\ref{fig:b}.
With this notation $$\calB'=\{t_{k_1,\ldots,k_s} \; | \; k_i\in\mathbb \Z \setminus \{0\},\; s\in\mathbb N,\;k_1\ge k_2\ge \cdots \ge k_s\}.$$

We will exhibit a similar basis of $\hsm(T)$, using an order relation similar to the one used in~\cite{DL}:
$$\calB''=\{t_{k_1, \ldots,k_s} \; | \; k_i\in\mathbb Z \setminus \{ 0\},s\in\mathbb N, 0>k_1\ge\ldots\ge k_l,0<k_{l+1}\le\ldots\le k_s\}.$$


We introduce an order relation on configurations of ovals with arrows, i.e. on diagrams with no crossings,
and such that there are no ovals with zero arrows.
This relation is defined in a lexicographical way by considering in this order:
\begin{enumerate}
	\item The number of all arrows counted as positive for all ovals independently of orientations.\label{e1}
\item The number of ovals.\label{e2}
\item A lexicographical ordering of ovals with positive arrows: from ovals with the smallest number of arrows to 
ovals with the largest number of arrows.\label{e3}
\item A lexicographical ordering of ovals with negative arrows in the same way as in 3), but taking absolute
values (i.e. the number of arrows) into account.\label{e4}
\end{enumerate}

To illustrate \ref{e3}, in $\calB$ one has, $t_1^2 t_3^2>t_1^2 t_2 t_4$, because $(1,1,3,3)>(1,1,2,4)$ (lexicographically).
To illustrate \ref{e4}, $t_{-2}^2 t_{-4}^2>t_{-2}^2 t_{-3} t_{-5}$ because $(|-2|,|-2|,|-4|,|-4|)>(|-2|,|-2|,|-3|,|-5|)$.

The order defined above becomes a total order when restricted to $\calB$ or $\calB''$. Indeed, if two
configurations of ovals have the same order, they will have the same series of arrows on ovals,
for example $(-3,-3,-1,1,4,4)$. Restricting to $\calB$ or $\calB''$ this determines completely the diagram.

Lemma~\ref{lemma:hsm_push} can be refined to take into account this ordering, in the following sense:

\begin{lemma}\label{lemma:push_order}
Let $D$ be a diagram with no crossings, with an oval $t_i$ containing no other ovals, next to
an oval $t_j$ (it can be nested in $t_j$ or not). Suppose that $0<|i|\le |j|$ if $i$ and $j$ have the
same sign, or $i<0<j$ otherwise. Then one can push $t_i$ through $t_j$ getting $v^{\pm 2i}$ times the diagram
in which the whole $t_i$ is pushed plus terms of lower order.
\begin{proof}

If $i<0<j$, then it follows from Lemma~\ref{lemma:hsm_push} that all terms will have less arrows than $D$,
except for the term corresponding to the whole $t_i$ being pushed, which comes with a factor $v^{-2i}$.
If $0<i\le j$, then, from the same lemma, it follows that we will have a term as before, this time with coefficient $v^{2i}$, 
plus terms for which some $a>0$ arrows will be transferred from $t_i$ to $t_j$. In the lexicographical order the change
will be from $(\ldots,i,\ldots,j,\ldots)$ to $(\ldots,i-a,\ldots,j+a,\ldots)$ which is of lower order.
Similarly for $j<i<0$ there will be a change from $(\ldots,j,\ldots,i,\ldots)$ to $(\ldots,j-a,\ldots,i+a,\ldots)$ which is again of lower order.
\end{proof}
\end{lemma}


\begin{theorem}\label{thm:new_basis_torus}
$\calB''$ is a basis of the free skein module $\hsm(T)$. 
\begin{proof}
Let $T=t_{k_1}\cdots t_{k_s}\in\calB$, $k_1\le k_2\le \cdots\le k_s$, $k_i\in \Z \setminus \{ 0\}$.
We construct first a function $F:\calB\to R\calB''$.
If $k_l<0$ and $k_{l+1}>0$, $F(T)$ will be equal to $t_{k_l,k_{l-1},\cdots,k_1,k_{l+1},\cdots,k_s}$ 
(which has the same order as $T$) times an invertible coefficient plus terms of lower order.
$F$ is defined by induction on the order. $F$ is the identity if the number of arrows is
$0$ (empty diagrams) or if the number of ovals is $1$.

If $k_s>0$ push all other ovals into $t_{k_s}$, using Lemma~\ref{lemma:push_order}.
We get $v$ to some power times $t_{k_s}$ around all other ovals plus terms of lower order.
Reexpress these terms of lower order in $\calB$ by pushing ovals out of $t_{k_s}$ again.
As we push through the oval with a maximum of arrows, 
Lemma~\ref{lemma:push_order} guarantees that the order cannot increase.
Thus $F$ of the terms that were reexpressed in $\calB$ will be of lower order than the
term with all ovals pushed into $t_{k_s}$ and the induction can be applied.

If $k_s<0$ push all ovals into $t_{k_1}$ to get again $v$ to some power times $t_{k_1}$
around all other ovals plus terms of lower order which, again, are reexpressed with terms
in $\calB$ of lower order.

Now repeat this process for ovals inside $k_s$ or $k_1$, pushing ovals into the oval
with maximum positive arrows (or maximum negative if there are no ovals with positive arrows).
For terms with lower order reexpress them in $\calB$: if some arrows were killed we will obviously
get terms with lower order also in $\calB$; if some arrows were transferred (by construction from ovals with less arrows
to ovals with more arrows) we will get terms with lower lexicographical order. 
It is clear that the order cannot increase back to the order of $T$, when pushing the ovals
which have lost some arrows out to get elements in $\calB$.

Continue until all ovals
are nested. We get at the end $v$ to some power times $t_{k_l,k_{l-1},\cdots,k_1,k_{l+1},\cdots,k_s}\in\calB''$
plus terms of lower order.

Now extend $F$ linearly from $R\calB$ to $R\calB''$. The matrix of $F$ with respect
to the ordered $\calB$ and $\calB''$ will be upper triangular with invertible elements
(powers of $v$) on the diagonal. This shows that $\calB''$ is a basis of $\hsm(T)$.
\end{proof}
\end{theorem}

Recall that the basis $\mathcal{B}_p$ of $\hsm(L(p,1))$ consists of diagrams with non-nested
ovals and the number of arrows $k_i$ on each of them satisfying $-\frac{p}{2}<k_i\le\frac{p}{2}$.

We want to exhibit a new basis of this skein module, $\calB_p''$, using a proof similar to that
of Theorem~\ref{thm:new_basis_torus}. Let: 
$$\calB_p''=\{t_{k_1,..,k_s} \; | \; k_i\in \Z \setminus \{0\},\; s\in\mathbb N,$$ $$0>k_1\ge\ldots\ge k_l>-\frac{p}{2}, \;0<k_{l+1}\le\ldots\le k_s\le \frac{p}{2}\}.$$

Thus $\calB_p''$ is $\calB''$ with arrows on ovals restricted to the interval $(-\frac{p}{2},\frac{p}{2}]$.
The order on $\calB''$ used in the proof of the preceding theorem restricts to an order on $\calB_p''$. We use it in the proof
of the next theorem.

\begin{theorem}\label{thm:new_basis_lp1}
$\calB_p''$ is a basis of the free skein module $\hsm(L(p,1))$.
\begin{proof}
We construct a function $F:\calB_p\to R\calB_p''$ in the same way as it was done in the proof of Theorem~\ref{thm:new_basis_torus},
having a similar property, namely that $F$ of an element $T$ in $\calB_p$ will be equal to the element of $\calB_p''$ consisting
of the same ovals as in $T$ but nested, times an invertible coefficient plus terms of lower order.
When reexpressing an element in $\calB_p$ it may happen that the number of arrows is reduced with $\Om_{(p)}$ moves but this lowers the order.

Extending $F$ linearly to $R\calB_p$, its matrix is again upper triangular with invertible elements on the diagonal,
from which it follows that $\calB_p''$ is a basis of $\hsm(L(p,1))$.
\end{proof}
\end{theorem}

\subsection*{Acknowledgements}
	The first author was supported by the Slovenian Research Agency grants J1-8131, J1-7025 and N1-0064.



\end{document}